\documentclass[reqno]{amsart}
\usepackage{amsmath, amsfonts, epsfig, amssymb}

{\theoremstyle{definition}
\newtheorem{definition}{Definition}[section]
\newtheorem{example}[definition]{Example}
}
{\theoremstyle{plain}%
  \newtheorem{theorem}{Theorem}
  \newtheorem{corollary}[definition]{Corollary}
  \newtheorem{proposition}[definition]{Proposition}
  \newtheorem{lemma}[definition]{Lemma}%
}
{\theoremstyle{remark}
\newtheorem{remark}[definition]{Remark}
}

\newcommand{\0}{\hat{0}}
\newcommand{\1}{\hat{1}}
\newcommand{\co}{\lessdot}
\newcommand{\coeq}{\leq$\raisebox{1.1pt}{$\negmedspace\!\cdot$}$\ }
\newcommand{\ga}{\gamma}
\newcommand{\de}{\delta}
\newcommand{\M}{M}
\newcommand{\compatible}{viable}
\newcommand{\compatibility}{viability}
\newcommand{\m}{m}
\newcommand{\ch}{c}
\newcommand{\R}{R_{\M}(\ch)}
\newcommand{\Q}{Q_{\M}(\m)}
\newcommand{\vez}{\vee^z}
\newcommand{\wedgy}{\wedge_y}
\newcommand{\wedgt}{\wedge_t}
\newcommand{\veu}{\vee^u}
\newcommand{\wedgu}{\wedge_u}
\newcommand{\nsn}{NS_n}
\newcommand{\bl}[1]{\langle #1 \rangle}

\begin{document}

\title{Poset Edge-Labellings and Left Modularity}

\author{Peter McNamara}
\address{Laboratoire de Combinatoire et d'Informatique Math\'ematique\\
Universit\'e du Qu\'ebec \`a Montreal\\
Case Postale 8888, succursale Centre-ville\\ 
Montr\'eal (Qu\'ebec) H3C 3P8\\
Canada}
\email{mcnamara@lacim.uqam.ca}

\author{Hugh Thomas}
\address{Fields Institute, 222 College St., Toronto, ON, M5T 3J1, Canada}
\email{hugh@math.unb.ca}


\begin{abstract}  
It is known that a graded lattice of rank $n$ is supersolvable if and only if 
it has an EL-labelling where the labels along any maximal chain are exactly
the numbers $1,2,\ldots,n$ without repetition.  These labellings are called
$S_n$ EL-labellings, and having such a labelling is also equivalent to 
possessing a maximal chain of left modular elements.  In the case of an
ungraded lattice, there is a natural extension of $S_n$ EL-labellings, 
called interpolating labellings. We show that admitting an 
interpolating labelling is again equivalent to possessing a maximal chain 
of left modular elements.
Furthermore, we work in the setting of an arbitrary bounded poset as all the
above results generalize to this case.  
We conclude by applying our results to show that 
the lattice of non-straddling partitions, 
which is not graded in general, has a maximal chain of left modular elements.
\end{abstract}

\maketitle


\section{Introduction}
An \emph{edge-labelling} of a poset $P$ is a map
from the edges of the Hasse diagram of $P$ to $\mathbb{Z}$.  
Our primary goal is to express certain classical properties of 
$P$ in terms of edge-labellings admitted by $P$.  
The idea of studying edge-labellings of posets goes back to \cite{St}.
An important milestone was \cite{Bj}, where 
A. Bj\"orner defined EL-labellings,
and showed that if a poset admits an EL-labelling, then it is 
shellable and hence Cohen-Macaulay.
We will be interested in a subclass of 
EL-labellings, known as $S_n$ EL-labellings.  
In \cite{St1}, R. Stanley
introduced supersolvable lattices and showed that they admit 
$S_n$ EL-labellings.  Examples of supersolvable lattices include
distributive lattices, the lattice of partitions of $[n]$, the 
lattice of non-crossing partitions of $[n]$ and the lattice of
subgroups of a supersolvable group (hence the terminology).
It was shown in \cite{Mc} that a finite graded lattice of rank $n$ is 
supersolvable if and only if it admits an $S_n$ EL-labelling.  In many 
ways, this characterization of lattice supersolvability in terms of
edge-labellings serves as the starting point for our investigations.

For basic definitions concerning partially
ordered sets, see \cite{ec1}.
We will say that a poset $P$ is \emph{bounded} if it contains a unique
minimal element and a unique maximal element, denoted $\0$ and $\1$
respectively.  
All the posets we will consider will be finite and bounded.
A chain of a poset $P$ is said to be \emph{maximal} if it is maximal under 
inclusion.
We say that $P$ is \emph{graded} if all the maximal chains of $P$ have the same
length, and we call this length the \emph{rank} of $P$.
We will write $x \lessdot y$ if $y$ covers $x$ in $P$ and
$x \coeq y$ if $y$ either covers or equals $x$.  
The edge-labelling $\ga$ of $P$ is
said to be an \emph{EL-labelling} if for any $y<z$ in $P$, 
\begin{enumerate}
\item[(i)] there is a unique unrefinable chain 
$y=w_0 \co w_1 \co \cdots \co w_r = z$ such that 
$\ga(w_0,w_1) \leq \ga(w_1,w_2) \leq \cdots \leq \ga(w_{r-1},w_r)$, and
\item[(ii)] the sequence of labels of this chain (referred to as the
\emph{increasing chain} from $y$ to $z$), when read from bottom 
to top, lexicographically precedes the labels of any other unrefinable chain
from $y$ to $z$.  
\end{enumerate}
This concept originates in \cite{Bj}; for the case where
$P$ is not graded, see \cite{BW1,BW2}.
If $P$ is graded of rank $n$ with an EL-labelling $\ga$, then $\ga$ is 
said to be an \emph{$S_n$ EL-labelling} if the labels along any maximal
chain of $P$ are all distinct and are elements of $[n]$.  In other words,
for every maximal chain 
$\0=w_0 \co w_1 \co \cdots \co w_n = \1$ of $P$, the map sending $i$ to 
$\ga(w_{i-1},w_i)$ is a permutation of $[n]$.
Note that the second condition in the definition of an EL-labelling is
redundant in this case.

\begin{example}\label{eg:distributive}
Any finite distributive lattice has an $S_n$ EL-labelling.
Let $L$ be a finite distributive lattice of rank $n$.  By
the Fundamental Theorem of Finite Distributive Lattices
\cite[p. 59, Thm. 3]{Bi}, that is
equivalent to saying that
$L = J(Q)$, the lattice of order ideals of some $n$-element poset $Q$.
Let $\omega : Q \to [n]$ be a linear extension of $Q$, i.e., any
bijection labelling the vertices of $Q$ that is
order-preserving (if $a<b$ in $Q$ then $\omega(a)<\omega(b)$).
This labelling of the vertices of $Q$ defines a labelling of the edges of
$J(Q)$ as follows.  If $y$ covers $x$ in $J(Q)$, then the order ideal
corresponding to $y$ is obtained from the order ideal corresponding to
$x$ by adding a single element, labeled by $i$, say.  Then we set
$\ga(x,y) = i$.  This gives us an $S_n$ EL-labelling for $L=J(Q)$.
Figure \ref{fig:distexample} shows a labelled poset and its lattice
of order ideals with the appropriate edge-labelling.
\begin{figure}
\center
\epsfxsize=60mm
\epsfbox{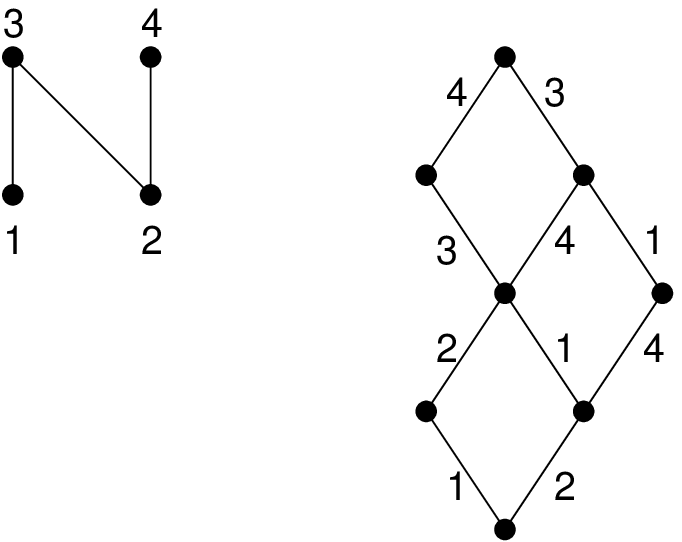}
\caption{}
\label{fig:distexample}
\end{figure}
\end{example}

A finite lattice $L$ is said to be 
\emph{supersolvable} if it contains a maximal chain,
called an 
\emph{M-chain} of $L$, which together with any other chain
in $L$ generates a distributive sublattice.  
We can label each such distributive sublattice by the method described in
Example \ref{eg:distributive} in such a way that the M-chain is the 
unique increasing maximal chain.
As shown in \cite{St1}, this will assign
a unique label to each edge of $L$ and the resulting global labelling
of $L$ is an $S_n$ EL-labelling. 

There is also a characterization of lattice supersolvability in
terms of
 left modularity.  Given an element $x$ of a finite lattice $L$, and
a pair of elements $y \leq z$, it is always true that
\begin{equation}\label{eq:modularineq}
(x \vee y) \wedge z \geq (x \wedge z) \vee y.
\end{equation}
The element $x$ is said to be \emph{left modular} if, for all $y \leq z$,
equality holds in \eqref{eq:modularineq}.  
Following A. Blass and B. Sagan \cite{BS}, we will say that a lattice itself is 
\emph{left modular} if it contains a left modular maximal chain, that is, a
maximal chain each of whose elements is left modular.  (One might
guess that we should define a lattice to be left modular if all of its
elements are left modular, but this is equivalent to the definition of 
a modular lattice.) 
As shown in \cite{St1}, 
any M-chain of a supersolvable lattice
is always a left modular maximal chain, and so supersolvable lattices are
left modular.  
Furthermore, it is shown by L. S.-C. Liu 
\cite{L} that if
$L$ is a finite graded lattice with a left modular maximal chain $\M$, then
$L$ has an $S_n$ EL-labelling with increasing maximal chain $\M$.
In turn, as shown in \cite{Mc}, this implies that $L$ is 
supersolvable, and so
we conclude the following.

\begin{theorem}\label{thm:gradedlattice}
Let $L$ be a finite graded lattice of rank $n$.  Then the following 
are equivalent:
\begin{enumerate}
\item $L$ has an $S_n$ EL-labelling,
\item $L$ is left modular,
\item $L$ is supersolvable.
\end{enumerate}
\end{theorem}

It is shown in \cite{St1} that if $L$ is upper-semimodular, then $L$ is 
left modular if and only if $L$ is supersolvable.  
Theorem \ref{thm:gradedlattice} is
a considerable strengthening of this.  Here we used $S_n$ EL-labellings
to connect left modularity and supersolvability. 
It is natural to ask for a more direct proof that (2) implies (3); such 
a proof has recently been 
provided by the second author in \cite{Th}.

Our goal is to generalize Theorem \ref{thm:gradedlattice} 
to the case when $L$ is not 
graded and, moreover, to the case when $L$ is not necessarily a lattice.
We now wish to define natural generalizations of $S_n$ EL-labellings and
of maximal left modular chains.

\begin{definition}
An EL-labelling $\ga$ of a poset
$P$ is said to be \emph{interpolating} if, for any 
$y \co u \co z$, either
\begin{enumerate} 
\item[(i)] $\ga(y,u) < \ga(u,z)$ or 
\item[(ii)] the increasing 
chain from $y$ to $z$, say $y = w_0 \co w_1 \co \cdots \co w_r = z$, has
the properties that its labels are strictly increasing and that
$\ga(w_0,w_1) = \ga(u,z)$ and $\ga(w_{r-1},w_r)=\ga(y,u)$.
\end{enumerate}
\end{definition}

\begin{example}
The reader is invited to
check that the labelling of the non-graded poset shown in Figure
\ref{fig:tamari} is an interpolating EL-labelling.
\begin{figure}
\center
\epsfxsize=50mm
\epsfbox{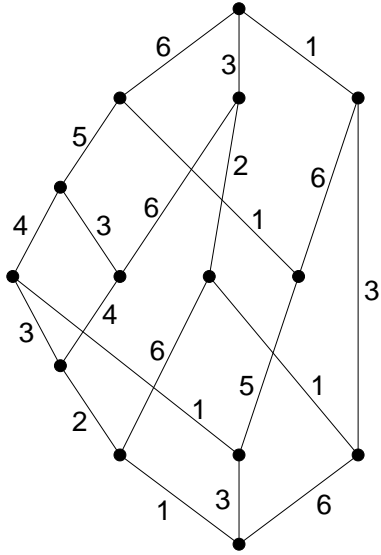}
\caption{The Tamari lattice $T_4$ and its interpolating EL-labelling}
\label{fig:tamari}
\end{figure}
In fact, the poset shown is the so-called ``Tamari lattice'' $T_4$.
For all positive integers $n$, there exists a Tamari lattice $T_n$
with $C_n$ elements, where $C_n = \frac{1}{n+1}\binom{2n}{n}$,
the $n$th Catalan number. 
More information on the Tamari lattice can be found in \cite[\S 9]{BW2},
\cite[\S 7]{BS}
and the references given there, and in
\cite[\S 3.2]{L}, where this interpolating EL-labelling appears.
The Tamari lattice is shown to have an EL-labelling in \cite{BW2} and
is shown to be left modular in \cite{BS}.
\end{example}

If $P$ is graded of rank $n$ and has an interpolating labelling $\ga$
in which the labels on the increasing maximal chain reading from bottom to top
are $1,2,\ldots n$, then we can check (cf. Lemma \ref{lem:labelsdistinct})
that $\ga$ is an
$S_n$ EL-labelling.  

Our next step is to define left modularity in the non-lattice case.
Let $x$ and $y$ be elements of $P$.  We know that $x$ and $y$ have at least
one common upper bound, namely $\1$.  If the set of common upper bounds 
of $x$ and $y$ has a least element, then we denote it by $x \vee y$.
Similarly, if $x$ and $y$ have a greatest common lower bound, then
we denote it by $x \wedge y$.

Now let $w$ and $z$ be elements of $P$ with $w, z \geq y$.  Consider the set
of common lower bounds for $w$ and $z$ that are also greater than or equal
to $y$.  Clearly, $y$ is in this set.  
If this set has a greatest element, then we denote it by $w \wedgy z$ and
we say that $w \wedgy z$ is well-defined (in $[y,\1]$).
We see that $(x \vee y)\wedgy z$ is well-defined
in the poset shown in Figure \ref{fig:snellable}, even though 
$(x \vee y) \wedge z$ is not.  Similarly, let $w$ and $y$ be 
elements of $P$ with $w, y \leq z$.  If the set
$\{u \in P\ |\ u \geq w,y \mbox{\ and\ } u \leq z\}$ has a least element, then
we denote it by $w \vez y$ and we say that $w \vez y$ is 
well-defined (in $[\0,z]$).  We will usually be interested in expressions
of the form $(x \vee y)\wedgy z$ and $(x \wedge z)\vez y$.  
The reader that is solely interested
in the lattice case can choose to ignore the subscripts 
and superscripts
on the meet and
join symbols.

\begin{definition}
An element $x$ of a poset 
$P$ is said to be \emph{\compatible} if, for all $y \leq z$ in $P$,
$(x \vee y) \wedgy z$ and $(x \wedge z) \vez y$
are well-defined.
A maximal chain of $P$ is said to be \compatible\ if each of its elements is
\compatible.
\end{definition}

\begin{example}
The poset shown in Figure \ref{fig:snellable} is certainly not a 
lattice but the reader can check that the increasing maximal chain 
is \compatible.
\begin{figure}
\center
\epsfxsize=55mm
\epsfbox{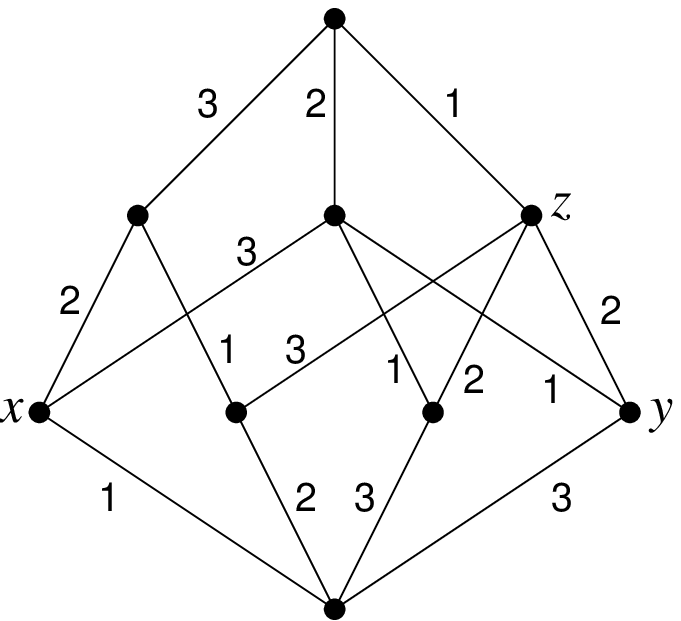}
\caption{}
\label{fig:snellable}
\end{figure}
\end{example}

\begin{definition} 
A \compatible\ element $x$ of a poset 
$P$ is said to be \emph{left modular} if, for all $y \leq z$ in $P$,
\[
(x \vee y) \wedgy z = (x \wedge z) \vez y.
\]
A maximal chain of $P$ is said to be left modular if each of its elements
is \compatible\ and left modular, and $P$ is said to be left modular if it possesses 
a left modular maximal chain.
\end{definition}

This brings us to the first of our main theorems.

\begin{theorem}\label{thm:theorem2}
Let $P$ be a bounded poset with a
left modular maximal chain $\M$. 
Then $P$ has an interpolating EL-labelling with $\M$ as its 
increasing maximal chain.
\end{theorem}

The proof of this theorem will be the content of the next section.  
In Section 3, we will prove the following converse result.

\begin{theorem}\label{thm:theorem3}
Let $P$ be a bounded poset with an interpolating EL-labelling.  The unique
increasing chain from $\0$ to $\1$ is a
left modular maximal chain.
\end{theorem}

These two theorems, when compared with Theorem \ref{thm:gradedlattice}, might 
lead one to ask about possible supersolvability results
for bounded posets that aren't graded lattices.  This problem is
discussed in Section 4.  
In the case of graded posets, we obtain a satisfactory result, namely
Theorem \ref{thm:posetsupersolvable}.  As a consequence, 
we have given an answer to the 
question of when a graded poset $P$ has an $S_n$ EL-labelling.  
This has ramifications
on the existence of a ``good 0-Hecke algebra action'' on the maximal
chains of the poset, as discussed in \cite{Mc}. 
However, it remains an open problem to appropriately extend the 
definition of supersolvability to ungraded posets.

An explicit 
application
of Theorem \ref{thm:theorem3} is the subject of Section 5.  As a variation on non-crossing
partitions and non-nesting partitions, we define non-straddling 
partitions.  Ordering the set of non-straddling partitions of $[n]$
by refinement gives a poset, denoted $\nsn$, 
that is generally a non-graded lattice.
We define an edge-labelling $\ga$ for $\nsn$ that is analogous
to the usual EL-labelling for the lattice of partitions of $[n]$.
In order to show that $\nsn$ is left modular, we then 
prove that $\ga$ is an interpolating EL-labelling.

\section{Proof of Theorem \ref{thm:theorem2}}

Throughout this section, we suppose that 
$P$ is a bounded poset with a left modular 
maximal chain $\M: \0 =x_0 \co x_1 \co \cdots \co x_n = \1$.
We want to show that $P$ has an interpolating EL-labelling.
Our approach will be as follows: we will begin by specifying
an edge-labelling $\ga$ for $P$ such that $M$ is an increasing chain
with respect to
$\gamma$.
We will then prove a 
series of lemmas which build on the \compatibility\ and
left modularity properties.  These
culminate with Proposition 
\ref{prop:labelsequal} which, roughly speaking, gives a more
local definition for $\ga$.  We will then be ready to show
that $\ga$ is an EL-labelling and is, furthermore, an 
interpolating EL-labelling.

We choose a label set 
$l_1 < \cdots < l_n$ of natural numbers.  (For most purposes, we can let
$l_i = i$.)  We define an edge-labelling $\ga$ on $P$ by setting
$\ga(y,z)=l_i$ for $y\co z$ if 
\[ (x_{i-1} \vee y)\wedgy z = y \mbox{\ \ and\ \ } 
(x_i \vee y) \wedgy z = z .\]
It is easy to see that $\ga$ is well-defined.  
We will refer to it as the labelling induced by $\M$ and the label set
$\{l_i\}$.
When $P$ is a lattice, this labelling appears, for example, in \cite{L,W}.
As in \cite{L}, we can give an equivalent definition of $\ga$ as follows.

\begin{lemma}
Suppose $y \co z$ in $P$.  Then $\ga(y,z)= l_i$ if and only if
\[ i = \min\{j\ |\ x_j \vee y \geq z \} = \max\{ j+1\ |\ x_j \wedge z \leq y\}.\]
\end{lemma}

\begin{proof}
That $i=\min\{j\ |\ x_j \vee y \geq z\}$ is immediate from the definition 
of $\ga$.  By left modularity, $\ga(y,z)=l_i$ if and only if
$(x_{i-1} \wedge z)\vez y = y$ and $(x_i \wedge z) \vez y = z.$
In other words, $x_{i-1}\wedge z \leq y$ and $x_i \wedge z \nleq y$.
It follows that $i = \max\{ j+1\ |\ x_j \wedge z \leq y \}$.
\end{proof}

\begin{lemma}\label{lem:repeatedjoin}
Suppose that $y \leq w \leq z$ in $P$ and let $x\in \M$.
Then 
$((x \wedge z)\vez y) \vez w$ is well-defined 
and 
equals $(x \wedge z) \vez w$.
Similarly, $((x \vee y)\wedgy z) \wedgy w$ is well-defined 
and 
equals $(x \vee y) \wedgy w$.
\end{lemma}

\begin{proof}
It is routine to check that, in $[\0,z]$, $(x \wedge z) \vez w$ 
is the least common upper bound for $w$
and $(x \wedge z)\vez y$, and that, in $[y,\1]$, 
$(x \vee y)\wedgy w$ is the greatest common lower bound lower bound for
$(x \vee y)\wedgy z$ and $w$.
\end{proof}

\begin{lemma}\label{lem:intervalmod}
Suppose that $t \leq u$ in $[y,z]$ and $x \in \M$.  
Let $w=(x \vee y) \wedgy z
=(x \wedge z) \vez y$ in $[y,z]$.  Then $(w \vez t) \wedgt u$ and
$(w \wedgy u) \veu t$ are well-defined elements of $[t,u]$ and are equal.
\end{lemma}

\begin{proof}
We see that, by Lemma \ref{lem:repeatedjoin},
\begin{equation*}
\begin{split}
(x \vee t) \wedgt u  &=  ((x \vee t) \wedgt z)\wedgt u 
= ((x\wedge z) \vez t)\wedgt u  \\ 
&=  (((x \wedge z)\vez y)\vez t)\wedgt u 
= (w \vez t) \wedgt u. \\
\end{split}
\end{equation*}
Similarly, 
\[ (x \wedge u)\veu t = (w \wedgy u)\veu t\]
But $(x \vee t)\wedgt u = (x \wedge u)\veu t$, yielding the result.
\end{proof}

\begin{lemma}\label{lem:covers}
Suppose $x$ and $w$ are \compatible\ and that $x$ is left modular in $P$. 
\begin{enumerate}
\item[(a)]
If $x \co w$ then for any $z$ in $P$ we have $x \wedge z \coeq w \wedge z$.
\item[(b)]
If $w \co x$ then for any $y$ in $P$ we have $w \vee y \coeq x \vee y$.
\end{enumerate}
\end{lemma}
Part (b) appears in the lattice case in \cite[Lemma 2.5.6]{L} and 
\cite[Lemma 5.3]{LS}.

\begin{proof}
We prove (a); (b) is similar. Assume, seeking a contradiction, that
$x \wedge z < u < w \wedge z$ for some $u \in P$. Now $u \leq z$ and
$u \leq w$.  It follows that $u \nleq x$.  

Now $x < x \vee u \leq w$.  Therefore, $w = x \vee u$.  So
\[ u = (x \wedge z) \vez u = (x \vee u)\wedgu z = w \wedge z ,\]
which is a contradiction.
\end{proof}

We now prove a slight extension of
\cite[Lemma 2.5.7]{L} and \cite[Lemma 5.4]{LS}.

\begin{lemma}\label{lem:wichain}
The elements of $[y,z]$ of the form $(x_i\vee y)\wedgy z$ form a
left modular maximal chain in $[y,z]$.
\end{lemma}

\begin{proof}
Lemma \ref{lem:intervalmod} gives the \compatibility\ and left modularity
properties. By Lemma \ref{lem:covers}(b), 
$x_i \vee y \coeq x_{i+1} \vee y$.  By Lemma \ref{lem:intervalmod} with
$z = \1$, we have that $x_i \vee y$ is left modular in $[y,\1]$.
Therefore, $(x_i \vee y) \wedgy z \coeq (x_{i+1} \vee y) \wedgy z$ by Lemma
\ref{lem:covers}(a).
\end{proof}

We are now ready for the last, and most important, of our preliminary
results.  
Let $[y,z]$ be an interval in $P$.
We call the maximal chain of $[y,z]$ from Lemma \ref{lem:wichain} the
\emph{induced}  
left modular maximal chain of $[y,z]$.
One way to get a second edge-labelling for $[y,z]$ would be to take
the labelling induced in $[y,z]$ by this induced maximal chain. 
We now prove that, for a suitable choice of label set, this
labelling coincides with $\ga$.

\begin{proposition}\label{prop:labelsequal}
Let $P$ be a bounded poset, 
$\0 =x_0 \co x_1 \co \cdots \co x_n = \1$ a
left modular maximal
chain and $\ga$ the corresponding edge-labelling with label set
$\{l_i\}$.  Let $y < z$, and define $c_i$ by saying
\begin{eqnarray*}
y & = & (x_{0} \vee y) \wedgy z = \cdots = (x_{c_1-1} \vee y) \wedgy z \\
  &   &\lessdot\ (x_{c_1} \vee y) \wedgy z = \cdots = (x_{c_2-1} \vee y) \wedgy z
  	\lessdot \cdots \\
 &  & \lessdot\ (x_{c_r} \vee y) \wedgy z = \cdots = (x_{n} \vee y) \wedgy z.
\end{eqnarray*}
Let $m_i = l_{c_i}$.  Let $\de$ be the labelling of $[y,z]$ induced
by its induced left modular maximal chain and the label set $\{m_i\}$.
Then $\de$ agrees with $\ga$ restricted to the edges of $[y,z]$.
\end{proposition}

\begin{proof}
Suppose $t \lessdot u$ in $[y,z]$.  
Using ideas from the proof of Lemma \ref{lem:intervalmod},
\begin{eqnarray*}
\de(t,u)=m_i & \Leftrightarrow & 
(((x_{c_i-1} \vee y)\wedgy z)\vez t)\wedgt u = t \mbox{\ and\ } \\
& & (((x_{c_i} \vee y)\wedgy z)\vez t)\wedgt u = u \\
& \Leftrightarrow & (x_{c_i-1} \vee t)\wedgt u = t 
\mbox{\ and\ } (x_{c_i} \vee
t)\wedgt u = u \\
& \Leftrightarrow & \ga(t,u)=l_{c_i}.
\end{eqnarray*}
\end{proof}

\begin{proof}[Proof of Theorem \ref{thm:theorem2}.]
We now know that the induced left modular chain in $[y,z]$ has (strictly)
increasing labels, say $m_1 < m_2 < \cdots < m_r$.
Our first step is 
to show that it is the 
only maximal chain with (weakly) increasing labels. 
Suppose that $y = w_0 \co w_1 \co \cdots \co w_r = z$ is the induced chain
and that $y = u_0 \co u_1 \co \cdots \co u_s = z$ is another chain with
increasing labels.  

If $s=1$ then $y \co z$ and the result is
clear.  Suppose $s \geq 2$.
By Proposition \ref{prop:labelsequal}, 
we may assume that the labelling on $[y,z]$ is induced by
the induced left modular chain $\{w_i\}$.
In particular, we have that $\ga(u_i, u_{i+1}) = m_l$ where 
$l = \min\{j\ |\ w_j \vez u_i \geq u_{i+1}\}$.  
Let $k$ be the least number such that $u_k \geq w_1$.  Then it is
clear that $\ga(u_{k-1},u_k)=m_1$.  Note that this is the smallest
label that can occur on any edge in $[y,z]$.
Since the labels on the chain $\{u_i\}$ are assumed to be increasing, 
we must have $\ga(u_0,u_1)=m_1$.  It follows that $w_1 \vez u_0 \geq u_1$
and since $y \lessdot w_1$, we must have $u_1 = w_1$.
Thus, by induction, the two chains coincide.
We conclude that the induced left modular maximal chain is the only
chain in $[y,z]$ with increasing labels.  

It also has the lexicographically
least set of labels.  To see this, 
suppose that $y = u_0 \co u_1 \co \cdots \co u_s = z$ is another chain 
in $[y,z]$.  We assume that $u_1 \neq w_1$ since, otherwise, we can just
restrict our attention to $[u_1, z]$.
We have 
$\ga(u_0, u_1) = m_l$, where 
$l=\min\{j\ |\ w_j \geq u_1\} \geq 2$ since $w_1 \ngeq u_1$.  Hence
$\ga(u_0,u_1) \geq m_2 > \ga(w_0, w_1)$.  
This gives that $\ga$ is an EL-labelling.
(That $\ga$ is an EL-labelling
was already shown in the lattice case in \cite{L,W}.)

Finally, we show that it is an interpolating EL-labelling.
If $y \co u \co z$ is not the induced left modular maximal chain in $[y,z]$, 
then let $y = w_0 \co w_1 \co \cdots \co w_r = z$ be the induced left
modular maximal chain.   
We have that $\ga(y,u) = m_l$ where 
\[
l=\min\{j\ |\ w_j\vez y \geq u\} =\min\{j\ |\ w_j \geq u\} = r
\]
since $u \co z$.
Therefore, $\ga(y,u) = m_r$.
Also, $\ga(u,z) = m_l$ where 
\[
l=\max\{j+1\ |\ w_j\wedgy z \leq u\} =\max\{j+1\ |\ w_j \leq u\} = 1
\]
since $y \co u$.
Therefore, $\ga(y,u) = m_1$, as required.
\end{proof}

\section{Proof of Theorem \ref{thm:theorem3}}

We suppose that $P$ 
is a bounded poset with an interpolating EL-labelling $\ga$.
Let $\0 = x_0 \co x_1 \co \cdots \co x_n = \1$ be the increasing chain
from $\0$ to $\1$ and let $l_i = \ga(x_{i-1},x_i)$.  
We will begin by
establishing some basic facts about interpolating labellings.
These results will enable us to show certain meets and joins exist by
looking at the labels that appear along particular increasing chains.
We will thus show that the $x_i$ are \compatible .
We will finish by showing that the $x_i$ are left 
modular, again by looking at the labels on increasing chains.

Let $y=w_0 \co w_1 \co \cdots \co w_r = z$.  Suppose that, for 
some $i$, we have
$\ga(w_{i-1},w_i) > \ga(w_i, w_{i+1})$.  Then the ``basic
replacement'' at $i$ takes the given chain and replaces the subchain
$w_{i-1} \co w_i \co w_{i+1}$ by the increasing chain from $w_{i-1}$ to
$w_{i+1}$.  The basic tool for dealing with interpolating labellings
is the following well-known fact about EL-labellings.

\begin{lemma}\label{lem:basicreps}
Let $y=w_0 \co w_1 \co \cdots \co w_r = z$.  Successively perform basic 
replacements on this chain, and stop when no more basic replacements
can be made.  This algorithm terminates, and yields the increasing chain
from $y$ to $z$.
\end{lemma}

\begin{proof}
At each step, the sequence of labels on the new chain lexicographically 
precedes the sequence on the old chain, so the process must terminate,
and it is clear that it terminates in an increasing chain.
\end{proof}
 
We now prove some simple consequences of this lemma.

\begin{lemma}\label{lem:labelsdistinct}
Let $\m$ be the chain $y=w_0 \co w_1 \co \cdots \co w_r = z$.  Then the 
labels on $\m$ all occur on the increasing chain from $y$ to $z$ and
are all different.  Furthermore, all the labels on the increasing 
chain from $y$ to $z$ are bounded between the lowest and highest labels on
$\m$.
\end{lemma}

\begin{proof}
That the labels on the given chain all occur on the increasing chain
follows immediately from Lemma \ref{lem:basicreps} and the fact that
after a basic replacement, the labels on the old chain all occur on 
the new chain.  Similar reasoning implies that the labels on 
the increasing chain are bounded between the lowest and highest labels
on $\m$.

That the labels are all different again follows from Lemma
\ref{lem:basicreps}.  Suppose otherwise.  By repeated basic replacements,
one obtains a chain which has two successive equal labels, which is not
permitted by the definition of an interpolating labelling.
\end{proof}

\begin{lemma}\label{lem:lessthanxi}
Let $z \in P$ such that there is some chain from $\0$ to $z$ all of
whose labels are in $\{l_1,\ldots,l_i\}$.  Then $z \leq x_i$.  Conversely,
if $z \leq x_i$, then all the labels on any chain from $\0$ to $z$ are
in $\{l_1,\ldots,l_i\}$.
\end{lemma}

\begin{proof}
We begin by proving the first statement.  By Lemma \ref{lem:labelsdistinct},
the labels on the increasing chain from $\0$ to $z$ are in
$\{l_1,\ldots,l_i\}$.  Find the increasing chain from $z$ to $\1$.
Let $w$ be the element in that chain such that all the labels below
it on the chain are in $\{l_1,\ldots,l_i\}$, and those
above it are in $\{l_{i+1},\ldots,l_n\}$.  Again, by Lemma
\ref{lem:labelsdistinct}, the increasing chain from $\0$ to $w$ has all
its labels in $\{l_1,\ldots,l_i\}$, and the increasing
chain from $w$ to $\1$ has all its labels in $\{l_{i+1},\ldots,l_n\}$.
Thus $w$ is on the increasing chain from $\0$ to $\1$, and so $w=x_i$.
But by construction $w \geq z$.  So $x_i \geq z$.

To prove the converse, observe that by Lemma \ref{lem:labelsdistinct}, no
label can occur more than once on any chain.  But since every label 
in $\{l_{i+1},\ldots,l_n\}$ occurs on the increasing chain
from $x_i$ to $\1$, no label from among that set can occur on any edge
below $x_i$.
\end{proof}

The obvious dual of Lemma \ref{lem:lessthanxi} is proved similarly:

\begin{corollary}
Let $z \in P$ such that there is some chain from $z$ to $\1$ all of
whose labels are in $\{l_{i+1},\ldots,l_n\}$.  Then $z \geq x_i$.  Conversely,
if $z \geq x_i$, then all the labels on any chain from $z$ to $\1$ are
in $\{l_{i+1},\ldots,l_n\}$.
\end{corollary}

We are now ready to prove the necessary \compatibility\ properties.

\begin{lemma}\label{lem:meetandjoin}
$x_i \vee z$ and $x_i \wedge z$ are well-defined for any $z \in P$
and for $i=1,2,\ldots,n$.
\end{lemma}

\begin{proof}
We will prove that $x_i \wedge z$ is well-defined.  The proof that
$x_i \vee z$ is well-defined is similar.
Let $w$ be the maximum element on the increasing chain from $\0$ to $z$
such that all labels on the increasing chain between $\0$ and $w$
are in $\{l_1,\ldots,l_i\}$.
Clearly $w \leq z$ and, by Lemma \ref{lem:lessthanxi}, $w \leq x_i$.

Suppose $y \leq z, x_i$.  It follows that all labels from $\0$ to $y$ are
in $\{l_1, \ldots,l_i\}$.  Consider the increasing chain from $y$ to
$z$.  There exists an element $u$ on this chain such that all the labels
on the increasing chain from $\0$ to $u$ are in $\{l_1, \ldots,l_i\}$
and all the labels on the increasing chain from $u$ to $z$ are in
$\{l_{i+1},\ldots,l_n\}$.  Therefore, $u$ is on the increasing chain
from $\0$ to $z$ and, in fact, $u = w$.  Also, we have that
$\0 \leq y \leq u = w \leq z$.  We conclude that $w$ is the greatest
common lower bound for $z$ and $x_i$.  
\end{proof}

\begin{lemma}
$\0 = x_0 \wedge z \leq x_1 \wedge z \leq \cdots \leq x_n \wedge z = z$,
after we delete repeated elements, is the increasing chain in $[\0,z]$.
Hence, $(x_i \wedge z) \vez y$ is well-defined for $y \leq z$.
Similarly, $(x_i \vee y)\wedgy z$ is well-defined.
\end{lemma}

\begin{proof}
From the previous proof, we know that $x_i \wedge z$ is the maximum 
element on the increasing chain from $\0$ to $z$ such that
all labels on the increasing chain between $\0$ and $x_i \wedge z$ are
in $\{l_1,\ldots,l_i\}$.  The first assertion follows easily from
this.  

Now apply Lemma \ref{lem:meetandjoin} to the bounded poset $[\0,z]$.
It has an obvious interpolating labelling induced from the interpolating
labelling of $P$.
Recall that our definition of the existence of $(x_i \wedge z) \vez y$
only requires it to be well-defined in $[\0,z]$.
The result follows.
\end{proof}

We conclude that the increasing maximal chain 
$\0 = x_0 \co x_1 \co \cdots \co x_n = \1$ of $P$ is \compatible.  
It remains
to show that it is left modular.

\begin{proof}[Proof of Theorem \ref{thm:theorem3}.]
Suppose that $x_i$ is not left modular for some $i$.  Then there exists
some pair $y \leq z$ such that $(x_i \vee y) \wedgy z > (x_i \wedge z)\vez y$.
Set $x=x_i$, $b=(x_i \wedge z)\vez y$ and $c=(x_i \vee y)\wedgy z$.
Observe that $d:=x\vee b \geq c$ while $a:= x \wedge c \leq b$. So the 
picture is as shown in Figure \ref{fig:proofmodular}.
\begin{figure}
\center
\epsfxsize=25mm
\epsfbox{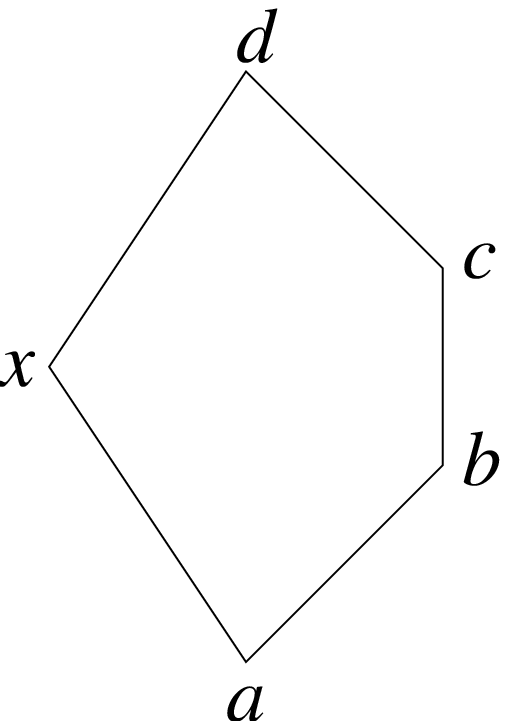}
\caption{}
\label{fig:proofmodular}
\end{figure}

By Lemma \ref{lem:lessthanxi}, the labels on the increasing chain from
$\0$ to $a$ are less than or equal to $l_i$.  Consider the increasing chain 
from $a$ to $c$.  Let $w$ be the first element along the chain.  If
$\ga(a,w) \leq l_i$, then by Lemma \ref{lem:lessthanxi}, $w \leq x_i$,
contradicting the fact that $a = x \wedge c$.  Thus the labels on the
increasing chain from $a$ to $c$ are all greater than $l_i$.  Dually, the
labels on the increasing chain from $b$ to $d$ are less than or equal to
$l_i$.  But now, by Lemma \ref{lem:labelsdistinct}, the labels on
the increasing chain from $b$ to $c$ must be contained in the labels on the
increasing chain from $a$ to $c$, and also from $b$ to $d$.  But there are
no such labels, implying a contradiction.  We conclude that the $x_i$ are
all left modular.
\end{proof}

We have shown that if $P$ is a bounded poset with an interpolating 
labelling $\ga$, then the unique increasing maximal chain $\M$ is a 
left modular maximal chain.  By Theorem \ref{thm:theorem2}, 
$\M$ then induces an
interpolating EL-labelling of $P$.  We now show that this labelling agrees
with $\ga$ for a suitable choice of label set, which is a special case of
the following proposition.

\begin{proposition}\label{prop:labellinguniqueness}
Let $\ga$ and $\de$ be two interpolating EL-labellings of a bounded poset $P$.
If $\ga$ and $\de$ agree on the $\ga$-increasing chain from $\0$ to $\1$,
then $\ga$ and $\de$ coincide.
\end{proposition}

\begin{proof}
Let $\m: \0 = w_0 \co w_1 \co \cdots \co w_r = \1$ be 
the maximal chain with the lexicographically first $\ga$ 
labelling among those chains for which $\ga$ and $\de$ disagree. 
Since $\m$ is not the 
$\ga$-increasing chain from $\0$ to $\1$, we
can find an $i$ such that $\ga(w_{i-1},w_i) > \ga(w_i,w_{i+1})$.  Let
$\m'$ be the result of the basic replacement at $i$ with respect to 
the labelling $\ga$.  Then the $\ga$-label sequence of $\m'$ 
lexicographically precedes that of $\m$, so $\ga$ and $\de$ agree on
$\m'$.  But using the fact that $\ga$ and $\de$ are interpolating,
it follows that they also agree on $\m$.  Thus they agree everywhere.
\end{proof}


\section{Generalizing Supersolvability}

Suppose $P$ is a bounded poset.  For now, we consider the case of $P$ being
graded of rank $n$.  We would like to define what it means for $P$ to 
be \emph{supersolvable}, thus generalizing Stanley's definition of
lattice supersolvability.  A definition of poset supersolvability with 
a different purpose appears
in \cite{W} but we would like a more general definition.  In particular,
we would like $P$ to be supersolvable if and only if 
$P$ has an $S_n$ EL-labelling.
For example, the poset shown in Figure \ref{fig:snellable}, while it doesn't
satisfy V. Welker's definition,  
should satisfy our definition.
We need to define, in the poset case, the equivalent of a sublattice 
generated by two chains.

Suppose $P$ has a \compatible\ maximal chain $\M$.  
Thus $(x \vee y) \wedgy z$ and $(x \wedge z)\vez y$ are well-defined 
for $x \in \M$ and $y \leq z$ in $P$.  Given any chain $\ch$ of $P$, we define 
$\R$ to be the smallest subposet of $P$ satisfying the following two
conditions:
\begin{enumerate}
\item[(i)] $\M$ and $\ch$ are contained in $\R$,
\item[(ii)] If $y \leq z$ in $P$ and $y$ and $z$ are in $\R$, then so are
	$(x \vee y) \wedgy z$ and $(x \wedge z)\vez y$ for any $x$ in $\M$.
\end{enumerate}

\begin{definition}
We say that a bounded poset $P$ is \emph{supersolvable} with M-chain $\M$
if $\M$ is a \compatible\ maximal chain and $\R$ is a distributive lattice 
for any chain $\ch$ of $P$.
\end{definition}

Since distributive lattices are graded, it is clear that a poset must
be graded in order to be supersolvable.  We now come to the main result of
this section.

\begin{theorem}\label{thm:posetsupersolvable}
Let $P$ be a bounded graded poset of rank $n$.  Then the following 
are equivalent:
\begin{enumerate}
\item $P$ has an $S_n$ EL-labelling,
\item $P$ is left modular,
\item $P$ is supersolvable.
\end{enumerate}
\end{theorem}

\begin{proof} Observe that for a graded poset, Lemma \ref{lem:labelsdistinct}
implies that an interpolating labelling is an $S_n$ EL-labelling, and
the converse is obvious.  Thus,
Theorems \ref{thm:theorem2} and \ref{thm:theorem3} 
restricted to the graded case give us that 
$(1) \Leftrightarrow (2)$.

Our next step is to show that $(1)$ and $(2)$ together imply $(3)$.  
Suppose $P$ 
is a bounded graded poset of rank $n$ with an $S_n$ EL-labelling.  
Let $\M$ denote the increasing maximal chain 
$\0=x_0 \co x_1 \co \cdots \co x_n = \1$   
of $P$.  We also know that $\M$ 
is \compatible\ and left modular and induces the same $S_n$ EL-labelling.  
Given any maximal chain $\m$ of $P$,
we define $\Q$ to be the closure of $\m$ in $P$ under basic replacements.
In other words, $\Q$ is the smallest subposet of $P$ 
which contains $\M$ and $\m$ and
which has the property that, if $y$ and $z$ are in $\Q$ with $y \leq z$, then
the increasing chain between $y$ and $z$ is also in $\Q$.
It is shown in 
\cite[Proof of Thm. 1]{Mc} 
that $\Q$ is a distributive lattice.
There $P$ is a lattice but the proof of distributivity doesn't use this
fact.  Now consider $\R$.  We will show that there exists a maximal 
chain $\m$ of $P$ such that $\R = \Q$.  
Let $\m$ be the maximal
chain of $P$ which contains $\ch$ and which has increasing labels between
successive elements of $\ch \, \cup \{\0,\1\}$.  
The only idea we need is that, for
$y \leq z$ in $P$, the increasing chain from $y$ to $z$ is given by
$y = (x_0\vee y) \wedgy z \leq (x_1\vee y)\wedgy z \leq \cdots \leq
(x_n \vee y)\wedgy z = z$, where we delete repeated elements.
This follows from Lemma \ref{lem:wichain} since the induced left modular 
chain in $[y,z]$ has increasing labels.  
It now follows that $\R = \Q$, and hence
$\R$ is a distributive lattice.

Finally, we will show that $(3) \Rightarrow (2)$.  We suppose that 
$P$ is a bounded supersolvable poset with M-chain $\M$.
Suppose $y \leq z$ in $P$ and let $\ch$ be the chain $y \leq z$.
For any $x$ in $\M$, $x \vee y$ is well-defined in $P$ 
(because $\M$ is assumed to be \compatible) and equals
the usual join of $x$ and $y$ in the lattice $\R$.  The same idea
applies to $x \wedge z$, $(x \vee y)\wedgy z$ and $(x\wedge z)\vez y$.
Since $\R$ is distributive, we have that 
\[
(x\vee y)\wedgy z 
= (x\vee y)\wedge z 
= (x\wedge z)\vee (y\wedge z) 
= (x \wedge z)\vee y
= (x \wedge z)\vez y
\]
in $\R$ and so $M$ is left modular in $P$.
\end{proof}

\begin{remark}
We know from Theorem \ref{thm:gradedlattice} that a graded lattice of
rank $n$ is supersolvable if and only if it has an $S_n$ EL-labelling.
Therefore, it follows from Theorem \ref{thm:posetsupersolvable} that
the definition of a supersolvable poset when restricted to graded
lattices yields the usual definition of a supersolvable lattice.  (Note 
that this is not {\it a priori} obvious from our definition of a supersolvable
poset.)
\end{remark}

\begin{remark}
The argument above for the equality of $\R$ and $\Q$ holds even if
$P$ is not graded.  However, in the ungraded case, it is certainly
not true that $\Q$ is distributive.  The search for a full generalization
of Theorem \ref{thm:gradedlattice} thus leads us to ask what can be
said about $\Q$ in the ungraded case.  Is it  a lattice?  Can we 
say anything even in the case that $P$ is a lattice?
\end{remark}

\section{Non-straddling partitions}

Let $\Pi_n$ denote the lattice of partitions of the set 
$[n]$ into blocks, where we order 
partitions by refinement: if $y$ and $z$ are partitions of $[n]$ 
we say that $y \leq z$
if every block of $y$ is contained in some block of $z$.  
Equivalently, $z$ covers $y$ in
$\Pi_n$ if $z$ is obtained from $y$ by merging two blocks
 of $y$.  Therefore, $\Pi_n$ is
graded of rank $n-1$.  $\Pi_n$ is shown to be supersolvable in 
\cite{St1} and hence has
an $S_{n-1}$ EL-labelling, which we denote be $\de$.  In fact, it will
simplify our discussion if we use the label set $\{2,\ldots ,n\}$ for
$\de$, rather than the label set $[n-1]$.  We choose
the M-chain, and hence the
increasing maximal chain for $\de$, to be the maximal chain 
consisting of the bottom element and those
partitions of $[n]$ whose only non-singleton block is $[i]$, 
where $2 \leq i \leq n$.  In the literature, $\de$ is often defined in
the following form, which can be shown to be equivalent.  If $z$ is obtained
from $y$ by merging the blocks $B$ and $B'$, then we set
\[
\de(y, z) = \max \{ \min B, \min B' \} .
\]
For any $x \in \Pi_n$,
we will say that $j \in \{2, \ldots n\}$ is a \emph{block minimum} in $x$
if $j = \min B$ for some block $B$ of $x$.  
In particular, we see that $\de(y,z)$ is the unique block minimum in $y$
that is not a block minimum in $z$.

Recall that a \emph{non-crossing} partition of $[n]$ is a partition with the
 property that if
some block $B$ contains $a$ and $c$ and some block $B'$ contains 
$b$ and $d$ with
$a < b < c < d$, then $B=B'$.  Again, we can order the set of non-crossing 
partitions of $[n]$
by refinement and we denote the resulting poset by $NC_n$.  This 
poset, which can be shown to be a lattice, has many nice properties and
has been studied extensively.  More information can be found in 
R. Simion's survey article \cite{Si} and the references given there.
Since $NC_n$ is a subposet of $\Pi_n$, we can consider $\de$ restricted to 
the edges of $NC_n$.  It was observed by Bj{\"{o}}rner and P. Edelman in 
\cite{Bj} that this gives an EL-labelling for $NC_n$ and we can easily see that this
EL-labelling is, in fact, an $S_{n-1}$ EL-labelling (once we subtract 1 from
every label).  

We are now ready to state our main definition for this section, which should be 
compared with the definition above of non-crossing partitions.

\begin{definition}
A partition of $[n]$ is said to be \emph{non-straddling} if whenever some
block $B$ contains $a$ and $d$ and some block $B'$ contains $b$ and $c$
with $a < b < c < d$, then $B=B'$. 
\end{definition}

This definition is also very similar to that of \emph{non-nesting} partitions, as
defined by A. Postnikov and discussed in \cite[Remark 2]{Re} and \cite{At}.
The only difference in the definition of non-nesting partitions is that we 
do not require $B=B'$ if there is also an element of $B$ between $b$ and $c$.
So, for example, $\{1,3,5\}\{2,4\}$ is a non-nesting partition in $\Pi_5$ but is not
a non-straddling partition.  We say that $\{1,3,5\}\{2,4\}$ is a 
\emph{straddling partition}, that $1 < 2 < 4 < 5$ is a \emph{straddle}, and that
the blocks $\{1,3,5\}$ and $\{2,4\}$ form a straddle. 

\begin{figure}
\center
\epsfxsize=135mm
\epsfbox{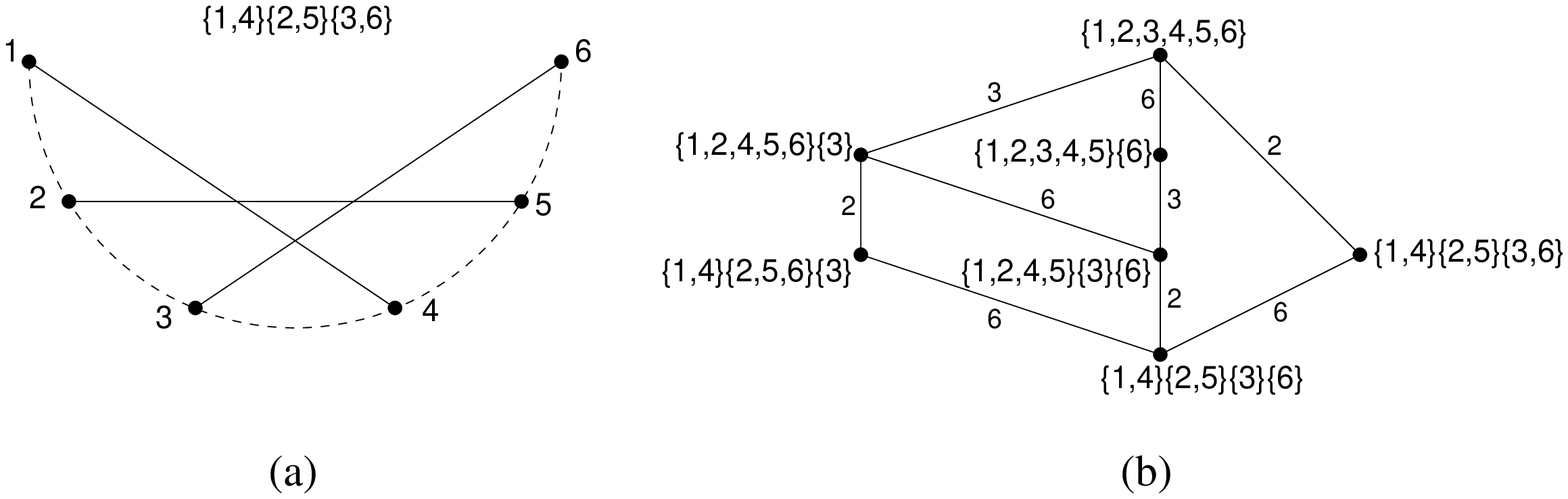}
\caption{}
\label{fig:nonstraddlinginterval}
\end{figure}

Let $\nsn$ be the subposet of $\Pi_n$ consisting of those partitions that are
non-straddling.  
To distinguish the interval $[x,y]$ in $\Pi_n$ from the interval $[x,y]$ in $\nsn$,
we will use the notation $[x,y]_{\Pi_n}$ and $[x,y]_{\nsn}$,
respectively.  We note that the meet in $\Pi_n$ of two non-straddling partitions
is again non-straddling, implying that $\nsn$ is a meet-semilattice. Since 
$\{1,2.\ldots, n\}$ is a top element for $\nsn$, we conclude that $\nsn$ is a 
lattice.  On the other hand, $\nsn$ is not graded.  For example, consider those
elements of $\Pi_6$ that cover $\{1,4\}\{2,5\}\{3,6\}$, as represented
in Figure \ref{fig:nonstraddlinginterval}(a).
$\{1,2,4,5\}\{3,6\}$, 
$\{1,3,4,6\}\{2,5\}$ and $\{1,4\}\{2,3,5,6\}$ are all straddling partitions, so
$\{1,4\}\{2,5\}\{3,6\}$ is covered in $NS_6$ by
$\{1,2,3,4,5,6\}$.
Figure \ref{fig:nonstraddlinginterval}(b) shows 
$[ \{1,4\}\{2,5\}\{3\}\{6\} ,\1]_{NS_6}$.  

Therefore, unlike $\Pi_n$ and $NC_n$, $\nsn$ cannot have an $S_{n-1}$ 
EL-labelling.  However, we can ask if it has an interpolating EL-labelling.
We see that the following three ways of defining an edge-labelling 
$\ga$ for $\nsn$ are
equivalent.  Observe that if $y \co z$ in $\nsn$, then $z$ is obtained from $y$
by merging the blocks $B_1, B_2, \ldots , B_r$ of $y$ into a single block $B$ in
$z$.  We set
\begin{eqnarray}
\ga(y,z) & = & \mbox{second smallest element of $\{\min B_1, \ldots , \min B_r\}$} \nonumber \\
& = & \mbox{smallest block minimum in $y$ that is not a block} \nonumber \\
&  & \mbox{minimum in $z$}  \nonumber \\
& = & \mbox{smallest edge label of $[y,z]_{\Pi_n}$  under the edge-labelling $\delta$.}  \label{eq:def3} 
\end{eqnarray}
See Figure \ref{fig:nonstraddlinginterval}(b) for examples.
Note that the label set for $\ga$ is $\{2, 3, \ldots , n\}$ and that if $r = 2$, then 
$\ga(y, z)$ equals $\de(y,z)$.  We see that the chain 
\[
\0 < \{1,2\}\{3\}\cdots \{n\} 
< \{1,2,3\}\{4\}\cdots \{n\}
< \cdots
< \{1,2,\ldots ,n-1\}\{n\} < \1
\]
is an increasing maximal chain in $\nsn$ under $\ga$.

\begin{theorem}\label{thm:gammainterpolating}
The edge-labelling $\ga$ is an interpolating EL-labelling for $\nsn$.
\end{theorem}

Applying Theorem \ref{thm:theorem3}, we get the following result:

\begin{corollary}
$\nsn$ is left modular.
\end{corollary}

In preparation for proving Theorem \ref{thm:gammainterpolating}, we wish to get a 
firmer grasp on $\nsn$.
Suppose $x, y \in \nsn$.  While the meet of $x$ and $y$ in $\nsn$ is just the meet of 
$x$ and $y$ in $\Pi_n$, the situation for joins is more complicated.  The next lemma, 
crucial to the proof that $\ga$ is an EL-labelling, helps us to understand important
types of joins.  From now on, unless otherwise specified, $x \vee y$ with $x, y \in \nsn$ will
denote the join of $x$ and $y$ in $\nsn$.  Furthermore, if 
$l_0 < l_1 < \cdots < l_r$ are block minima in $y$, then $\bl{l_i}$ will denote the block of
$y$ with minimum element $l_i$, and $\bl{l_0} \cup \bl{l_1} \cup \cdots \cup \bl{l_r}$ will denote
the minimum element $z \in \nsn$ for which the elements of $\bl{l_0}, \bl{l_1}, \cdots , \bl{l_r}$ are
all in a single block.  Note that $z$ is well-defined, since it is the meet of all those elements of
$\nsn$ that have the required elements in a single block. 

\begin{lemma}\label{lem:merging}
Suppose $l_0 < l_1 < \cdots < l_r$ are block minima in $y$ and that 
\[
 y \vee (\bl{l_0} \cup \bl{l_1})
 = y \vee (\bl{l_0} \cup \bl{l_1} \cup \cdots \cup \bl{l_r}).
\]
Then 
\[
y \vee (\bl{l_i} \cup \bl{l_j})
 = y \vee (\bl{l_0} \cup \bl{l_1} \cup \cdots \cup \bl{l_r}).
\]
for any $0 \leq i < j \leq r$.
\end{lemma}

In words, this says that if merging the blocks $\bl{l_0}$ and $\bl{l_1}$ in $y$ requires us to 
merge all of $\bl{l_0}, \bl{l_1}, \ldots , \bl{l_r}$, then merging any two of these blocks 
also requires us to merge all of them.  

\begin{proof}
The proof is by induction on $r$, with the result being trivially true when $r=1$.
While elementary, the details are a little intricate.  To gain a better understanding,
 the reader may wish to treat the proof as an exercise.
If $i < j < r-1$, then by the induction assumption and the hypothesis that
$y \vee (\bl{l_0} \cup \bl{l_1})
 = y \vee (\bl{l_0} \cup \bl{l_1} \cup \cdots \cup \bl{l_r}$, we have
\[
y \vee (\bl{l_i} \cup \bl{l_j}) = y \vee (\bl{l_0} \cup \bl{l_1} \cup \cdots \cup \bl{l_{r-1}})
= y \vee (\bl{l_0} \cup \bl{l_1} \cup \cdots \cup \bl{l_r}),
\]
as required.  
Therefore, it suffices to let $j = r$.  

Since 
$y \vee (\bl{l_0} \cup \bl{l_1}) = y \vee (\bl{l_0} \cup \bl{l_1} \cup \cdots \cup \bl{l_{r-1}})
= y \vee (\bl{l_0} \cup \bl{l_1} \cup \cdots \cup \bl{l_r})$, we know that 
$\bl{l_0} \cup \bl{l_1} \cup \cdots \cup \bl{l_{r-1}}$ forms a straddle with $\bl{l_r}$.
There are two ways in which this might happen.

Suppose we have $a < b < c < d$ with $a, d \in \bl{l_0} \cup \bl{l_1} \cup \cdots \cup \bl{l_{r-1}}$
and $b, c \in \bl{l_r}$.  Suppose $d \in \bl{l_s}$ in $y$.  Then, since $l_s < l_r \leq b < c$, we
have that $l_s < b < c < d$ is a straddle in $y$, which contradicts $y \in \nsn$.

Secondly, suppose we have $a < b < c < d$ with $a, d \in \bl{l_r}$
and $b, c \in \bl{l_0} \cup \bl{l_1} \cup \cdots \cup \bl{l_{r-1}}$.
Suppose  $b \in \bl{l_s}$ and $c \in \bl{l_t}$.  Now $c > b > a \geq l_r > l_s, l_t$.
If $s=t$ then $y$ has a straddle, so we can assume that $l_s \neq l_t$ and that
$l_i \neq l_t$, with the argument being similar if $l_i \neq l_s$.
If $l_i < l_t$, then $l_i < l_t < c < d$ is a straddle when we merge blocks 
$\bl{l_i}$ and $\bl{l_r}$ in $y$.
Therefore,
\begin{equation}\label{eq:merging}
y \vee (\bl{l_i} \cup \bl{l_r}) = y \vee (\bl{l_i} \cup \bl{l_t} \cup \bl{l_r})
= y \vee (\bl{l_0} \cup \bl{l_1} \cup \cdots \cup \bl{l_r})
\end{equation}
by the induction assumption.
If $l_i > l_t$, then $l_t < l_i < l_r < c$ is a straddle when we merge blocks
$\bl{l_i}$ and $\bl{l_r}$ in $y$, also implying \eqref{eq:merging}.
\end{proof}

\begin{lemma}\label{lem:firstx}
Suppose $y < z$ in $\nsn$ and that $[y, z]_{\Pi_n}$ has edge labels
$l_1 < l_2 < \cdots < l_s$ under the edge-labelling $\de$.  
\begin{enumerate}
\item[(i)] There is exactly one edge of the form $y \co w$ 
with $\ga(y, w) = l_1$ in $[y, z]_{\nsn}$.  
\item[(ii)]
On any unrefinable chain $y \co u_0 \co u_1 \co \cdots \co u_k = z$ in $\nsn$, the 
label $l_1$ has to appear.
\end{enumerate} 

\end{lemma}

\begin{proof}
(i) We first prove the existence of $w$.
Let $l_0$ be the minimum of the block of $z$ containing $l_1$
and set $w = y \vee (\bl{l_0} \cup \bl{l_1})$.  Suppose $y < u \leq w$.  We know
$w$ is obtained from $y$ by merging the blocks $\bl{l_0}, \bl{l_1}, \bl{l_{i_1}},
\bl{l_{i_2}}, \bl{l_{i_r}}$, for some $0 \leq r < s$.  Applying Lemma \ref{lem:merging},
we get that $u = w$ and so $y \co w$.  By definition of $\ga$, we have that 
$\ga(y, w) = l_1$.  

It remains to prove uniqueness.  Suppose $w' \in \nsn$ with $y \co w'$ in $[y,z]$.
If $\ga(y, w') = l_1$, then we see that the blocks $\bl{l_0}$ and $\bl{l_1}$ must 
be merged in $w'$.  Therefore, these two blocks are merged in $w \wedge w'$, which
is thus greater than $y$.
Since $y \co w, w'$, we conclude that $w = w'$.

(ii) Consider the chain $y=u_0 < u_1 < \cdots < u_k=z$ as a chain in $\Pi_n$.
Since $\de$ is an $S_{n-1}$ EL-labelling for $\Pi_n$ (once we subtract 1 from every
label), the label $l_1$ has to appear on every maximal chain of $[y,z]_{\Pi_n}$.
It particular, it has to appear in one of the intervals $[u_i, u_{i+1}]_{\Pi_n}$ 
for $0 \leq i < k$.
Therefore, by \eqref{eq:def3}, we get that $\ga(u_i, u_{i+1}) = l_1$ for some
$0 \leq i < k$.  
\end{proof}

\begin{proposition}
The edge-labelling $\ga$ is an EL-labelling for $\nsn$.
\end{proposition}

\begin{proof}
Consider $y, z \in \nsn$ with $y < z$.  Suppose $[y,z]_{\Pi_n}$ has
edge labels $l_1 < l_2 < \cdots < l_s$.  By \eqref{eq:def3}, these are the only edge
labels that can appear in $[y,z]_{\nsn}$.  
We now describe a recursive construction of an unrefinable chain 
$\lambda: y = w_0 \co w_1 \co \cdots \co w_k = z$ in $\nsn$.
We let $w_1$ be the $w$ of Lemma \ref{lem:firstx}, i.e. $w_1$ is that unique element 
of the interval $[y,z]$ in $\nsn$ that covers $y$ and satisfies $\ga(y, w_1)=l_1$. 
Obviously, the labels in the interval $[w_1, z]$ are all greater than $l_1$.  
Now we apply the same argument in the interval $[w_1, z]$ to define $w_2$ and repeat until 
we have constructed all of $\lambda$.  Clearly, $\lambda$ is then an increasing chain.
By Lemma \ref{lem:firstx}(i), it has the lexicographically least set of labels.  By Lemma
\ref{lem:firstx}(ii), it is the only increasing chain from $y$ to $z$.
\end{proof}

\begin{proof}[Proof of Theorem \ref{thm:gammainterpolating}.]
Suppose we have $y \co u \co z$ in $\nsn$ with $\ga(y, u) > \ga(u, z)$.  Let 
$y = w_0 \co w_1 \co \cdots \co w_k = z$ be the unique increasing chain of $[y, z]$ in 
$\nsn$.  By Lemma \ref{lem:firstx}, we know that 
$\ga(w_0, w_1)=\ga(u,z)=l_1$, the smallest edge label
of $[y, z]_{\Pi_n}$.  

To show that $\ga(y, u) = \ga(w_{k-1}, w_k)$, we have to work considerably harder.  
We will continue to write $\bl{m}$ to denote the block of $y$ whose minimum is 
$m$ and we suppose that $u$ is obtained from $y$ by merging blocks
$\bl{m_0}, \bl{m_1}, \ldots , \bl{m_s}$ of $y$, with
$m_0 < m_1 < \cdots < m_s$.  We will write $\bl{l}_u$ to denote 
the block of $u$ whose minimum is $l$, and we suppose that $z$ is obtained
from $u$ by merging blocks $\bl{l_0}_u, \bl{l_1}_u, \ldots , \bl{l_r}_u$, 
with $l_0 < l_1 < \cdots < l_r$.
With the structure of the chain $y \co u \co z$ thus fixed, we now can deduce
information about the structure of the increasing chain.

If $l_0$ and $m_0$ are distinct and are both block minima in $z$, then all the 
$l_i$'s and $m_j$'s are distinct.  It follows that 
$z$ is obtained from $y$ by merging blocks $\bl{l_0}, \bl{l_1}, \ldots , \bl{l_r}$ and 
separately merging blocks $\bl{m_0}, \bl{m_1}, \ldots , \bl{m_s}$.  
Since $y \co u \co z$, we get that $k=2$ and 
$\ga(y, u) = \ga(w_1, z) = m_1$. 
We assume, therefore, that $m_0 = l_i$ for some $0 \leq i \leq r$.

As usual, we let $w_1 = y \vee (\bl{l_0} \cup \bl{l_1})$.  Now consider 
\[
w = y \vee (\bigcup_{i: \, l_i < m_1} \bl{l_i}) .
\]
Since $l_1 < m_1$, we know that $w \geq w_1$.  Let $B$ denote the 
the block of $w$ containing all $\bl{l_i}$ 
satisfying $l_i < m_1$.  Since $m_0 < m_1$, we know that 
$m_0 \in B$.  In fact, if we can show that $m_1 \not\in B$, then
we can now complete the proof.  Indeed, assume $m_1 \not\in
B$ and let $w' = w \vee (B \cup \bl{m_1})$. 
Now $w'$ has 
$m_0$ and $m_1$ in the same block and so satisfies $w' \geq u$, since
$u = y \vee (\bl{m_0} \cup \bl{m_1})$.  
Also, $w'$ has $l_0$ and $l_1$ in the same block and so satisfies
$w' \geq z$, since $z = u \vee (\bl{l_0} \cup \bl{l_1})$. 
Hence, $w' = z$.  By Lemma \ref{lem:merging} (substitute $w$ for $y$, 
and $m_1 < \cdots < m_s$ for $l_1 < \cdots < l_r$), we see that $w \co w'$.
Now $\ga(w, w') = m_1$, while
the edge labels of $[y, w]_{\Pi_n}$ all come from
the set $\{l_i\ | \ l_i < m_1\}$, implying that $w$ is on the increasing chain
between $y$ and $z$.  Therefore, $w = w_{k-1}$ and so $\ga(w_{k-1},w_k)= 
\ga(y,u)$.  

It remains to show that $m_1 \not\in B$.  
In fact, we will show that $m_j \not\in B$ for any $j \geq 1$.  
Consider the set:
\[
\tilde{B} = \bigcup_{i: \, l_i < m_1} \bl{l_i}
\]
We will show that $\tilde B$ does not form a straddle with any 
$\bl{m_j}$ for $j\geq 1$. From that, it follows immediately that 
$B=\tilde B$, and therefore that $m_1 \not\in B$, as desired.  

For $j\geq 1$, if $\bl{m_j}$ is a singleton, then $\tilde B$ does not
form a straddle with $\bl{m_j}$.  So suppose that $|\bl{m_j}|\geq 2$.
Let $m_j'$ denote
the second smallest element of $\bl{m_j}$.  Observe the following:
\begin{itemize}
\item If $\bl{m_0}$ contains an element greater than $m_j'$, then $\bl{m_0}$ and $\bl{m_j}$
form a straddle in $y$, which is impossible.
\item If $\bl{m_0}$ has more than one element between $m_j$ and $m_j'$, then we can draw
the same conclusion.
\item Consider those $l_i < m_1$ with $l_i \neq m_0$.  If $\bl{l_i}$ 
contains an element greater than $m_j$, then 
$\bl{l_i}_u$ forms a straddle in $u$ with $\bl{m_0} \cup \bl{m_1} \cup \cdots \cup \bl{m_s}$, 
which is impossible. 
\end{itemize}
Combining these three observations, we see that $\tilde B$ contains
no elements greater than $m_j'$, and at most one element 
between $m_j$ and $m_j'$.
In particular, it does not form a straddle with $\bl{m_j}$, as desired. 
\end{proof}


\section*{Acknowledgements}
The authors would like to thank Andreas Blass, Bruce Sagan, Richard Stanley, 
Volkmar Welker, and the anonymous referees for helpful comments.


\end{document}